\documentclass[12pt,a4paper,twoside]{amsart}

\usepackage[T1]{fontenc}
\input xypic
\usepackage{amsmath}
\usepackage{amssymb, graphics, color}
\usepackage{multicol}
\usepackage{mathrsfs}

\newcommand{\vc}[2]{\begin{pmatrix} #1 \\ #2 \end{pmatrix}}

\newcommand{\vol}[0]{\operatorname{Vol}}

\def\newtexttt{}
\def\newtextt{}
\def\newtext {}

\def\0{{\bf {0}}}

\makeatletter
\renewcommand{\subsection}{\@startsection{subsection}{2}{0mm}{\baselineskip}%
{-\fontdimen2\font plus -\fontdimen3\font minus \fontdimen4\font}%
{\normalfont\normalsize\bfseries}}
\renewcommand{\@seccntformat}[1]{\csname the#1\endcsname .\hspace{0.5em}}
\makeatother

\usepackage{amsmath}
\newtheorem {lemma}{Lemma}

\newtheorem {theorem}[lemma]{Theorem}
\newtheorem {definition}[lemma]{Definition}

\newtheorem {proposition}[lemma]{Proposition}

\def\M{{M}}

\def\tilde{\widetilde}

\def \mbeq {\begin {eqnarray}}
\def \meeq {\end {eqnarray}}

\def \bfo {\begin {displaymath} }
\def \efo {\end {displaymath} }

\def \beq {\begin {eqnarray}}
\def \eeq {\end {eqnarray}}
\def \ba {\begin {eqnarray*}}
\def \ea {\end {eqnarray*}}

\def \R {{\Bbb {R}}}

\def \H2s {H^{s+1}_0(\partial \M\times [0,T/2])}

\def \supp {\hbox{supp }}
\def \diam {\hbox{diam }}

\def \det {\hbox{det}}
\def\bra{\langle}
\def\cet{\rangle}

\def \e {\epsilon}

\def \a {\alpha}

\def \pa0 {\partial _0}

\def \p {\partial}

\def \tilde{\widetilde}

\def \e {\epsilon}

\def \A {{\mathcal A}}
\def \F {{\mathcal F}}

\title[Focusing waves in unknown media]{Focusing waves in unknown media by modified time reversal iteration}

\author[Dahl]{Matias Dahl}
\address{
Matias Dahl and Matti Lassas, Institute of Mathematics, P.O.Box
1100, 02015 Helsinki University of Technology, Finland }

\author[Kirpichnikova]{Anna Kirpichnikova}
\address{
Anna Kirpichnikova, School of Mathematics, The University of
Edinburgh, JCMB Mayfield Road Edinburgh EH9 3JZ, UK }

\author[Lassas]{Matti Lassas}

\date{\today}
\begin{document}
\maketitle
\begin{abstract}
  We study the wave equation in a bounded domain or on a compact
  Riemannian manifold with boundary.  Assume that we are given the
  hyperbolic Neumann-to-Dirichlet map on the boundary corresponding to
  physical boundary measurements.  We consider how to focus waves,
  that is, how to find Neumann boundary values so that at a given time
  the corresponding {\newtexttt wave converges to a delta distribution
    $\delta_y$ while the time derivative of the wave converges to
    zero. Such boundary value are generated by an iterative sequence
    of measurements.} In each iteration step we apply time reversal
  and other simple operators to measured data and compute boundary
  values for the next iteration step.  {\newtexttt The key feature of
    the algorithm is that it does not require knowledge of the
    coefficients in the wave equation, that is, the material
    parameters inside the media. However, we assume that the point $y$
    where the wave focuses is known in travel time coordinates.}
\end{abstract}


\section{Introduction}
Let us consider the wave equation in a bounded domain $M$,
\begin{eqnarray}
\label{eq: Wave}
\begin{cases}
  u_{tt}(x,t)+\A u(x,t)=0,\quad \hbox{ in } M\times \R_+,\\
  u|_{t=0}=0,\quad u_t|_{t=0}=0,  \\
  \p_\nu u|_{\p M\times \R_+}=f,
\end{cases}
\end{eqnarray} 
where $A$ is a 2nd order elliptic partial differential
operator. 

In this paper we show how to construct Neumann boundary values $f$
such that at time $T$, the wave $(u(T),u_t(T))$ is arbitrarily close
to $(c\delta_{y},0)$, where $\delta_{y}$ is the Dirac delta
distribution at a chosen point $y\in M$. {\newtextt We call such waves
  \emph{focusing waves}.}  To find such boundary values, we only
assume that we can make physical measurements from the boundary of
$M$. For given Neumann boundary values we can measure the Dirichlet
boundary values of the wave.  {\newtext A focusing wave can then be
  generated by an iterative sequence of measurements. In each
  iteration step we apply time reversal and other simple operators to
  measured data and compute boundary values for the next iteration
  step.

The iteration algorithm in this paper is closely related to time
reversal methods.  Let us therefore shortly discuss the underlying
idea and the {\newtextt usually used} approximations behind these methods.  As a simple
example, let us consider a domain $M$ in $\R^3$, and suppose that we
can measure waves and generate sources on the boundary of $M$.
Let us first assume that there is a theoretical point source at $y\in M$, 
and we measure the wave and/or its normal derivative at the boundary
of $M$. Assume further that we record this signal, reverse it in
time, and re-emitted into $M$, see \cite{FinkMain}.
Then one can show (assuming certain approximations hold, see
\cite{FinkMain,Fink2006,FinkCassereau}), that the re-emitted wave will
travel like the original wave, but as if time were running backwards.
This causes the re-emitted wave to focus near $y$.
}  

{\newtext
  This principle can also be used for imaging. To find a small
  scatterer $D$ in a {\newtextt relatively} homogeneous domain $M$, one sends a wave
  into $M$.  If the scatterer is small and the single scattering
  approximation is justified, 
  the scattered wave corresponds to a wave produced by a point source
  at $D$.  If we record this scattered signal at the boundary,
  reversed it in time, and re-emit it into the domain, it will focus
  at the scatterer.  Furthermore, this focusing has been observed to
  be quite stable under perturbations of the medium.  Thus, if the
  re-emitted wave is simulated (by computational means) in homogeneous
  media, it will focus at the location of $D$.  In this way a small
  scatterer can be found using relatively simple computational
  methods.  }
The above basic idea has been refined in various ways. If the target
area contains multiple scatterers, an iteration scheme can be used to
focus the wave on any of the scatterers \cite{PF}.

{\newtext Besides imaging, time reversal can be used to focus a wave
  onto a scatterer, say, inside the human body. One application of
  this is \emph{litotripsy}, where one breaks down a kidney or bladder
  stone using a focusing ultrasonic wave.  Another application is
  \emph{hyperthermia}, where a cancer is destroyed by an excessive
  heat dose generated by a focusing wave.  Let us point out that for
  the wave equation, there are various methods to estimate material
  parameters in travel time coordinates from boundary measurements.
  These methods are, however, quite unstable
  \cite{AKKLT,KatsudaKurylevLassas}. Therefore they might not be suitable
  for hyperthermia, where safety is crucial. 
An important question is therefore how to focus waves in unknown media.
}

For reviews and extensions on time reversal, see seminal papers of
M.\ Fink, \cite{FinkD,Fink1992,Fink2006}.
Time reversal methods have been intensively studied in random
heterogeneous media where the statistic of the random media is known,
see e.g. \cite{Bal2,Bal3,Bal4,Papa1,Papa3}.  For time
reversal in chaotic cavities, see \cite{room1}.  For related
analysis on time reversal methods, see also
\cite{Bardos,Bardos2,CIL,Kliba,PTF}.

Let us describe the key features of the algorithm in this paper.
First, to focus a wave we do not require knowledge of the material
parameters inside the media. We only assume that the coordinates of
the point are known in travel time coordinates. This means that
focusing can be done in the same coordinates in which imaging is done.
Thus, as the algorithm for focusing does not require media parameters
obtained from imaging, {\newtext 
  errors in imaging do not accumulate into errors in focusing.  }
Second, the algorithm can focus waves near an area having no
scatterers.  Third, the algorithm is computationally cheap. In a
sense, all computations are done in the media; there is no need to
solve the wave equation, cf \cite{I}.  {\newtexttt We will assume that
  the medium is linear, non-dispersive, non-dissipative, and
  frequency-independent. However, we do not need any other
  approximations like the single scattering approximations to prove
  that the algorithm works.  }

The limitations of the present algorithm is that we  assume
selfadjointness of operator $\A$ and that time $T$ is large
enough.  Moreover, the point $y$ where the wave focuses need to be
specified in  travel time coordinates  unless operator $\A$ is known.

The present work is a continuation of \cite{BKLS} where a similar
iterative scheme was introduced, for which $u(T)$ focuses to a
delta-distribution, but the time derivative $u_t(T)$ is uncontrolled.
The present work can also be seen as also a generalization of
so-called retrofocusing in control theory, where the aim is to produce
boundary sources giving the same final state as a boundary sources
sent before in the medium, see \cite{Jon-DeHoop,KL_edin}.  {\newtext
  The methodology in this paper arises from boundary control methods
  used to study inverse problems in hyperbolic equations
  \cite{B1,B2,BK,KK,KKL,KKLima,KKLM}.  }

{\newtext
The outline of this work is as follows.  In Section \ref{Terminology}
we introduce notation and review some relevant results from control theory.
We also define the boundary operators that
are needed in the iteration scheme.  In Section \ref{mainResults} we
describe the main results (Theorems \ref{Main F} and \ref{Cor 1}) and
outline their proofs, and in Section \ref{proofSec} we prove these
results.
}

\section{Definitions}
\label{Terminology}
{\newtext
We assume that $M\subset\R^m$ $(m\ge 1)$ is the closure of an open
{\newtextt $C^\infty$-smooth bounded set with non-empty smooth boundary $\p M$
or a $C^\infty$-smooth compact manifold with boundary. Furthermore, we
assume that $M$ is equipped with  a $C^\infty$-smooth} Riemannian metric
$g=\sum_{jk} g_{jk} \,dx^j \otimes dx^k$. Elements of the inverse
matrix of $g_{ij}$ are denoted by $g^{ij}$. 
Let $\mathrm{d}V_g$ be the {\newtexttt smooth measure}
\begin{eqnarray*}
  \mathrm{d}V_g &=& \vert g(x)\vert^{1/2}dx^1\wedge\cdots \wedge dx^m,
\end{eqnarray*}
where $|g|=\det([g_{jk}])$. Then $L^2(M)$ is defined by the inner product
\ba 
  \bra u,v\cet = \int_M u(x)v(x)\, \mathrm{d}V, 
\ea 
where $\mathrm{d}V = \mu \mathrm{d}V_g$ and 
$\mu \in C^\infty(M)$ is a fixed strictly positive function on $M$.

In wave equation \eqref{eq: Wave}, we assume $\A$ represents the 
most general formally selfadjoint elliptic partial differential operator 
with respect to the above inner product \cite{KKL}. 
In local coordinates, $\A$ has the form 
\ba
  \label{A} 
\A v=- \sum_{j,k=1}^m \frac{1}{\mu(x) |g(x)|^{1/2}} \frac {\p}{\p
x^j}\left( \mu(x) |g(x)|^{1/2}g^{jk}(x)\frac {\p v}{\p x^j} \right)+q(x)v, 
\ea 
where $q$ is a smooth function $q\colon M\to \R$. 
For example, if $\mu=1$ and $q=0$
then $\A$ reduces to the  Riemannian Laplace operator.
Let us point out that $\A$ represents media that is linear, 
{\newtexttt non-dissipative}, non-dispersive, and frequency-independent. 
}

On the boundary, operator $\p_\nu$ is defined by
\ba 
\p_\nu v=\sum_{j=1}^m \mu(x)\nu^j\frac {\p}{\p x^j}v(x)
\ea 
where $\nu(x)=(\nu^1,\nu^2,\dots,\nu^m)$ is the unit interior normal vector
satisfying
$\sum_{j,k=1}^m g_{jk}\nu^j\nu^k=1$. 
{\newtext
To integrate functions on $\p M$ we 
use the {\newtexttt measure} $\mathrm{d}S$ on $\p M$ induced by $\mathrm{d}V_g$. 
}
If $B\subset\p\M\times \mathbb{R}_+$, we define
$$
  L^2(B)=\{f\in L^2(\p\M\times\mathbb{R}_+) : \operatorname{supp}(f)\subset
B\}
$$ 
identifying functions and their zero continuations. 

With these assumptions, the wave equation has a solution whenever
$f\in L^2(\p M\times \R_+)$, and we denote this solution by $u^f$.
The map $f\mapsto u^f$ is linear over $\mathbb{R}$, and $\p_t u^f =
u^{\p_t f}$ when $f,\p_t f \in L^2(\p M, \mathbb{R}_+)$.

Let $d(x,y)$ be the \emph{geodesic distance} corresponding to $g$.
The metric $d$ is also called the \emph{travel time metric} because it
describes how solutions to equation \eqref{eq: Wave} propagate.  
By the  finite velocity
of wave propagation, (see \cite{HormanderIV}) we have that if
$\Gamma\subset \p M$ is open, and $f\in
L^2(\Gamma\times\mathbb{R}_+)$, then at time $t>0$, solution $u^f$ is
supported in the \emph{domain of influence}
\begin{eqnarray*}
\label{domain_of_influence} 
  M(\Gamma,t)=\{x\in M\ : \hbox{d}(x,\Gamma)\leq t\}.
\end{eqnarray*}
{\newtext
The \emph{diameter} of $M$ is defined as
\begin{eqnarray*}
\label{diameter_definition}
  \operatorname{diam}(M) = \max{\{\hbox{d}(x,y) : x,y\in M\}}.
\end{eqnarray*}
The \emph{characteristic function} of a set $S$ is denoted by $\chi_S$.
}

\subsection{Controllability for wave equation} 
\label{Sec: Controllability}
The seminal result implying controllability is Tataru's unique
continuation result \cite{Ta1,Ta3}.

\begin{proposition} [Tataru] 
\label {th:3.4} 
Let $u \in H^1_{\operatorname{loc}}(M\times \R_+)$ 
be a solution of the wave equation 
$$ 
  u_{tt}(x,t)+\A u(x,t)=0.
$$ 
Assume that 
\begin{eqnarray*}
  u|_{\Gamma\times (0,2\tau)}=0,\quad \p_\nu u|_{\Gamma \times (0,2\tau)}=0,
\end{eqnarray*}
where  $\Gamma\subset \p M$ is an non-empty open set and $\tau>0$.
Then 
\bfo 
  u(x,\tau)=0,\ \p_tu(x,\tau)=0 \hbox{ for }x\in M(\Gamma,\tau). 
\efo 
\end{proposition}

Using Tataru's unique continuation result, one can prove
the following controllability results {\newtextt (The proof is  postponed 
to section \ref{proofSec})}:

\begin{proposition}[Approximate global controllability]
\label{lemma_Control 1} 
If $T> \operatorname{diam}(\M)$, then the linear subspace
$$
  \left\{(u^f (T), u^f_t(T)) : \,f\in C_0^\infty
(\p\M\times\mathbb{R}_+) \right\}
$$
is dense in $H^1(M)\times L^2(\M)$.
\end{proposition}

This result yields the following controllability result, see e.g.
\cite{KKL} and references therein.

\begin{proposition}[Approximative local controllability] 
{\newtext 
\label{lemma_Control_2} 
Let $\tau>0$, 
let $\Gamma_1, \ldots, \Gamma_J \subset \p M$ be open non-empty sets, and
let $0<s_k<\tau$,  $k=1,\ldots, J$. Suppose
\begin{eqnarray*}
  B=\bigcup^{J}_{j=1} \Gamma_j\times (\tau-s_j,\tau),
\end{eqnarray*}
{\newtexttt
and $P$ is multiplication by the characteristic function $\chi_B$,
}
\begin{eqnarray}
\label{Pdefeq}
  P\colon L^2(\p\M\times(0,\tau)) &\rightarrow& L^2(\p\M\times(0,\tau)), \\
\nonumber
  f(x,t) &\mapsto & \chi_{B}(x,t)\,f(x,t).
\end{eqnarray}
Then the linear subspace
$$
  \left\{u^{Ph}(\tau) : h \in L^2 (\p M, (0,\tau) ) \right\}
$$
is dense in $L^2(N)$, where $N=\bigcup^{J}_{j=1}M(\Gamma_j,s_j)$.
}
\end{proposition}

\subsection{Operators for boundary sources}
In this section we introduce operators for manipulating 
boundary sources. These will be needed both in the proof of 
the main result and in the iteration scheme. 

For initial boundary value problem \eqref{eq: Wave} we
  define the \emph{non-stationary Neumann-to-Dirichlet map}
  (or \emph{response operator}) $\Lambda$ by setting 
\begin{eqnarray*}
  \Lambda f&=&u^f|_{\p M\times \R_+}, \quad f\in L^2(\p M\times \R_+).
\end{eqnarray*}
In other words, we solve the wave equation \eqref{eq: Wave} for a
boundary source $f$, and measure boundary values for the solution 
$u^f$ when $t> 0$. 
In this work, we only need the finite time Neumann-to-Dirichlet map, 
\begin{eqnarray*}
  \Lambda_{2T} f=u^f|_{\p\M\times(0,2T)},
\end{eqnarray*}
where $T>0$. By \cite{tataruREG} the map 
$$
  \Lambda_{2T}\colon L^2(\p M\times(0,2T))\rightarrow H^{1/3}(\p M \times
(0,2T))
$$ 
is bounded, where $H^s(\p M\times(0,2T))$ is the Sobolev space 
on $\p M\times(0,2T)$. 

For $f\in L^2(\p M, (0,2T))$, let 
\begin{eqnarray*}
 R_{2T}f(x,t)&=&f(x,2T-t), \\\nonumber
J_{2T}f(x,t)&=& \int_{[0,2T]} J_{2T}(s,t)f(x,s)ds,
\end{eqnarray*}
where $J_{2T}(s,t)=\frac 12\chi_{L}(s,t)$ and 
\begin{eqnarray*}
  L&=&\{ (s,t)\in \R_+\times\R_+:\ t+s\leq 2T,\ s>t \}. 
\end{eqnarray*}
We call $R_{2T}$ the \emph{time reversal map}, and $J_{2T}$ the 
\emph{time filter map} \cite{BKLS}. 
{\newtext
On $L^2(\p \M\times [0,2T])$ with the measure $dS(x)dt$,
the adjoint of $\Lambda_{2T}$ is \cite{BKLS}, 
\ba 
  \Lambda^*_{2T}=R_{2T}\Lambda_{2T} R_{2T}.
\ea  
}
For $f\in L^2(\p M, [0,2T])$, let  
$$
  Q_{2T} f=\int_0^{2T} g(t,s)f(x,s)ds,$$ 
be {\newtexttt the time filter operator,
where $g\colon (0,2T)^2 \to \R$, 
$$
g(t,s)=\frac{1}{2(e^{4T}-1)}
  \begin{cases}
    (e^{-t}+e^{t})(e^{4T}e^{-s}+e^{s}),\quad t<s,\\
    (e^{-s}+e^s)(e^{4T}e^{-t}+e^t),\quad t>s,
  \end{cases}\, \, 
$$ 
is the Green's function for the problem
$$
  \begin{cases}
    (1-\p_t^2)g(t,s)=\delta(t-s),\\ 
    \p_t g|_{t=0}=0,\quad \p_t g|_{t=2T}=0, \quad s\in (0,2T).
  \end{cases}
$$
}
Next we consider $\Lambda_{2T},\,R_{2T},\,J_{2T},\, Q_{2T}$ as
operators such that 
$$
  \Lambda_{2T},\,\,\,R_{2T},\,\,\,J_{2T}, \,\,\,Q_{2T} \colon L^2(\p\M \times[0,2T])
\rightarrow L^2(\p\M\times[0,2T]).
$$ 
Below, we often  denote $R_{2T},J_{2T},$ and $Q_{2T}$
by $R,J,$ and $Q$.
For $f,h\in L^2(\p \M\times [0,2T])$ the 
\emph{Blagovestchenskii identity}
states that
\beq
\label {4.60b}
\nonumber 
\int_M u^f(T)u^h(T)\,dV =\int_{\p M\times
[0,2T]} (Kf)(x,t) h(x,t)\, dS(x) dt, 
\eeq 
where $K\colon L^2\to L^2$ is the bounded operator 
\ba
K=K_{2T}=R_{2T}\Lambda_{2T}R_{2T}J_{2T}-J_{2T}\Lambda_{2T}.  
\ea 

For a proof, see e.g. \cite{BKLS}.  The importance of this identity is
that it shows that the inner product of {\newtexttt solutions $u^f(T)$
  and $u^h(T)$ can be calculated from the boundary.  Namely, on the
  right hand side of the Blagovestchenskii identity, $dS$ is the
  Riemannian surface volume on $\p\M$, and $K$ is defined in terms of
  the Neumann-to-Dirichlet map $\Lambda_{2T}$ and simple operators on
  boundary values like time reversal.  The intrinsic Riemannian
  surface volume $dS$ on $\p M$ is determined by $\Lambda_{2T}$.
  Namely, by Tataru's unique continuation principle, the Schwartz
  kernel of $\Lambda_{2T}$ is supported in
$$
  E=\{(x,t,x',t')\in(\p M\times[0,2T])^2:\,t-t' \ge d(x,x')\},
$$
and the boundary $\p E$ is in the support. The set $\p E$ determines
the distances of points $z,z'\in\p M$ with respect to the intrinsic
metric of the boundary $(\p M,g_{\p M})$.
}

\section{Iterations and main results} 
\label{mainResults}

\subsection{Cutoff of wave} 
In this section we describe Theorem \ref{Main F} which can be seen as a
lemma used in the proof of 
Theorem \ref{Cor 1}.

Let $X$ be the Banach space
$$
  X=L^2(\p M\times [0,2T])\times Y, \quad
  Y = H^1((0,2T); L^2 (\p M)),
$$
such that the inner product on $X$ is
\ba
  \left\bra \vc{h_1}{a_1}, \vc{h_2}{a_2} \right\cet_X &=& \bra h_1, h_2\cet_{L^2}
+ \bra a_1, a_2\cet_{L^2}
+ \bra \p_t a_1, \p_t a_2\cet_{L^2}.
\ea

\begin{definition}
\label{cutDef}
Let $T>2 \operatorname{diam}(M)$, and let
\begin{eqnarray*}
  B=\bigcup^{J}_{j=1} \Gamma_j\times (T-s_j,T),
\end{eqnarray*}
where $\Gamma_1, \ldots, \Gamma_J \subset \p M$ are open non-empty
sets, and $0<s_k<T$, $k=1,\ldots, J$. Let $P=\chi_B$ be the
multiplication with the characteristic function of $B$ defined as in
equation \eqref{Pdefeq}, and let $L\colon X\to X$ be the operator
\begin{eqnarray}
\label{Edef}
L = 
\begin{pmatrix}
  1 & 0 \\
  0 & Q \\
\end{pmatrix} 
 \begin{pmatrix}
  2PKP & -PK \\
  -KP & K-\p_t K\p_t \\
\end{pmatrix}.
\end{eqnarray}
Let $\alpha\in(0,1)$, and 
let $\omega>0$ be such that $2(1+\Vert L \Vert_X)<\omega$, and let
\begin{eqnarray*}
  S=(1-\frac {\alpha}{\omega})I-\frac {1}{\omega}L.
\end{eqnarray*}
If $f\in L^2(\p\M\times \mathbb{R}_+)$ be a boundary source, 
we define a sequence
$\vc{h_n}{a_n}=\vc{h_j(\alpha)}{a_j(\alpha)}\in X$, $n=1,2,\ldots$ by
\begin{eqnarray}
\label{firstItScheme}
\begin{cases}
  \vc{h_0}{a_0} = \displaystyle \frac{1}{\omega}\vc{PKf}{0}, \\
  \vc{h_{n}}{a_{n}} = \vc{h_0}{a_0} + S\vc{h_{n-1}}{a_{n-1}}, \quad n=1,2,\ldots.
\end{cases}
\end{eqnarray}
\end{definition}

\begin{theorem}[Cutoff of wave]
\label{Main F} 
Let $a_1(\alpha), a_2(\alpha),\ldots$,  be as 
in Definition \ref{cutDef}. Then the sequence converges in $Y$, 
$$
  \lim_{n\to \infty} a_{n}(\alpha)  = a(\alpha),
$$
and function $a(\alpha) \in Y$ on the right hand side satisfies
\ba
  \lim_{\alpha\to 0} \vc{u^{a(\alpha)}(T)}{u^{a(\alpha)}_t(T)} = 
\vc{\chi_{N} u^f(T)}{0}, 
\ea
where both  limits are in $L^2(M)$ and $N$ is the domain of influence
$$
    N = \bigcup_{k=1}^J M(\Gamma_k, s_k).
$$
\end{theorem}
Note that here $\omega$ may depend on $\alpha$. For instance,
we can choose $\omega=1/\alpha$.

Let us emphasize that the novelty of this theorem is the explicit
iteration scheme for $a(\alpha)$ depending only on boundary
measurements.  The scheme depends on operators $J,Q,P$, and $K$ that
can be calculated from the boundary of $M$.  The first three are
simple operators like integration and restriction.  Operator
$K=R_{2T}\Lambda_{2T}R_{2T}J_{2T}-J_{2T}\Lambda_{2T}$ involve time
reversal $R_{2T}$, time filtering $J_{2T}$ and two evaluations of the
Neumann-to-Dirichlet map $\Lambda_{2T}$ which corresponds to
two physical measurements. Hence, the first order
approximation of $a(\alpha)$ requires $2$ physical measurements. After
that, each additional term requires $10$ additional measurements.
Thus, for finite approximation of $a(\alpha)$ we only need finitely
many evaluations of the Neumann-to-Dirichlet map.



The full proof of Theorem \ref{Main F} is given in Section \ref{Sec: ProofOfFirstTheorem}. Let us here outline the main ideas.
For $\alpha\in(0,1)$, boundary sources $h(\alpha), a(\a)$ are defined as the 
minimum of the functional 
\beq 
\nonumber
\mathcal F(h,a,\a) &=& \|u^f(T)-u^{Ph}(T)\|^2_{L^2(M)} \\
\nonumber
        & & \, + \|u^{Ph} (T) - u^a(T)\|^2_{L^2(M)}  + \|u^a_t(T)\|^2_{L^2(M)} \\
\label{F intro definition}
        & & \, + \alpha(\|h\|^2_{L^2(\p M\times [0,T])}+\|a\|^2_{L^2(\p M\times [0,T])}+\|\p_ta\|^2_{L^2(\p M\times [0,T])}), 
\eeq 
In the sequel, when there is no danger of misunderstanding, we denote
the $L^2$-norms in the spaces $L^2(M)$, $L^ 2(\p M\times [0,T])$ etc.\
just by $\|\,\cdotp\|$.  In Lemma \ref{lemma_7}, we use convexity to
prove that for each $\alpha$, there is a unique minimum $h(\alpha),
a(\a)$, and by studying the Fr\'echet derivative of $\F$ we find a
linear equation (equation \eqref{Eq: main}) for this minimum.
{\newtexttt
In Lemma \ref{lemmaXYZ} we show that iteration scheme
\eqref{firstItScheme} converges to a von Neumann sum that represents
the solution to equation \eqref{Eq: main}. }
That minimizer $a(\a)$ satisfies the sought limit is proven in Lemma
\ref{finalLemma}. The key step in the proof is to use the
approximative controllability results from Section \ref{Sec:
  Controllability} to show that the first terms in $\F$ can be
arbitrarily close to $\|(1-\chi_N)u^f(T)\|^2$ and the next two terms
can be made arbitrarily small.


\subsection{Focusing of wave} 
To understand how one can focus waves using Theorem
  \ref{Main F} suppose we have sets $B\subset \tilde B\subset \p M\times
  [0,T]$ (defined in terms of $\Gamma_i$ and $s_i$ as in Definition 
  \ref{cutDef}). 
  Then Theorem \ref{Main F} implies that there are boundary 
  sources $a(\alpha)$ and $\tilde a(\alpha)$ such that
\begin{eqnarray*}
  \lim_{\alpha\to 0}
  \vc{  u^{a(\alpha)}(T)}{u^{a(\alpha)}_t(T)} &=& \vc{\chi_{N} u^f(T)}{0}, \\
  \lim_{\alpha\to 0}
  \vc{  u^{\tilde a(\alpha)}(T)}{u^{\tilde a(\alpha)}_t(T)} &=& \vc{\chi_{\tilde N} u^f(T)}{0},
\end{eqnarray*}
where the domains of influences satisfy $N\subset \tilde N$ at time $T$. 
As the solution operator $f\mapsto u^f$ is linear and commutes with $\p_t$,
solution $b(\alpha) = \tilde a(\alpha) - a(\alpha)$ satisfies
\begin{eqnarray*}
  \lim_{\alpha\to 0}
  \vc{  u^{b(\alpha)}(T)}{u^{b(\alpha)}_t(T)} &=& \vc{\chi_{\tilde N\setminus N} u^f(T)}{0}.
\end{eqnarray*}
That is, in the limit, {\newtexttt the} solution corresponding to $b$
is supported in $\tilde N \setminus N$. In the proof we construct $P$
and $\tilde P$ such that $\tilde N \setminus N$ is a family of sets
that shrink onto a chosen point $\widehat x$.  By further scaling $b$
with a suitable constant depending on the volume of $\tilde N
\setminus N$, we obtain the delta distribution.

To formulate the result, let us introduce some notation. 
By $\gamma_{x,\xi}$ we mean a geodesic in $(M,g)$ parametrized
by the arclength such that  $\gamma_{x,\xi}(0)=x$,
$\dot\gamma_{x,\xi}(0)=\xi$, and $\|\xi\|_g=1$.  Let $\nu=\nu(z)$,
$z\in \p M$ be the interior unit normal vector to $\p M$.  
{\newtexttt Then} there is a
critical value $\tau(z)>0$, such that for $t<\tau(z)$ the geodesic
$\gamma_{z,\nu}([0,t])$ is the unique shortest geodesic from its
endpoint $\gamma_{z,\nu}(t)$ to $\p M$, and for $t>\tau(z)$ it is no
longer a shortest geodesic. 
{\newtexttt We will not consider the degenerate case $t=\tau(z)$.}

\begin{definition} 
\label{focusDef}
Let $T> 2\operatorname{diam}(M)$, let $\widehat{x}=\gamma_{\widehat
  z,\nu}(\widehat T)$, where $\widehat z\in \p M$, and $0<\widehat T<T$.
Let $\Gamma_j\subset \p M$ for $j=1,2,\dots$ be open sets around $z$,
such that $\Gamma_j \supset \overline \Gamma_{j+1}$ and
$\bigcap_{j=1}^\infty\Gamma_j=\{\widehat z\}$.

Suppose $f\in C^\infty_0(\p M\times \R_+)$. 
Let $a_n(\alpha, \varepsilon) \in Y$ be functions 
obtained from the iteration in Definition \ref{cutDef} when $B$ is the set
$$
  B(\varepsilon) = \p M \times \left(T-(\widehat{T}-\varepsilon), T\right),
$$
$\alpha\in (0,1)$, and $\varepsilon>0$ is sufficiently small. 
Similarly, let $a_n(\alpha, j, \varepsilon) \in Y$ be functions 
obtained from the iteration in Definition \ref{cutDef} when $B$ is the set
$$
  B(j,\varepsilon) = \left( \p M \times \left(T-(\widehat{T} - \varepsilon),T\right) \right) \,\cup\, \left(\Gamma_j \times \left(T-(\widehat{T}+\varepsilon),T \right)\right),
$$
$\alpha\in (0,1)$, $j=1,2,\ldots$ and $\varepsilon>0$ is 
sufficiently small. 

Under these assumptions, let 
$$
  b_n(\alpha, j, \varepsilon) =\varepsilon^{-\frac{m+1}{2}} \left(
a_n(\alpha, j, \varepsilon) - a_n(\alpha, \varepsilon)
\right) \in Y.
$$
\end{definition}

Theorem \ref{Cor 1} is the main result of this paper.

\begin{theorem}[Focusing wave]
\label{Cor 1} Let $\widehat z\in \p M$, $\widehat T$,
$b_n(\alpha, j,\varepsilon)$ for $n=1,2, \ldots$ be as in Definition
\ref{focusDef}. Then functions $b_j$ converge in $Y$, 
$$
  \lim_{n\to \infty} b_{n}(\alpha, j, \varepsilon)  = b(\alpha, j, \varepsilon).
$$
Moreover, if 
$\widehat T<\tau(\widehat z)$
then functions $b(\alpha, j, \varepsilon) \in Y$ satisfy
\beq\label{eq: A0} 
  \lim_{\varepsilon\to 0^+} \lim_{j\to \infty} \lim_{\alpha\to 0^+}
\vc{  u^{b(\alpha,j,\varepsilon)}(T)}{u^{b(\alpha,j,\varepsilon)}_t(T)} = 
 C(\widehat{x}) u^f(T,\widehat{x}) \vc{ \delta_{\widehat{x}}}{0},
\eeq
where the inner two limits are in $L^2(M)$ and the outer limit is in
$\mathscr{D}'(M)$, and $C(\widehat x)$ is defined in equation \eqref{Cdefeqn}.

{\newtextt
If
$\widehat T>\tau(\widehat z)$, the limit (\ref{eq: A0}) is zero. }
\end{theorem}

Let us make three comments about this theorem. First, the condition
$\widehat T<\tau(\widehat z)$ means that $\widehat z\in \p M$ is the
closest boundary point to $\widehat x$ and
$d(\widehat{x},z)=d(\widehat{x},\p M)$. If a point $\hat x$ has a
unique closest boundary point, we say that it is admissible. For
example, on the closed disc, all points are admissible except the
center. For a general $(M,g)$ the set of points that are not
admissible has measure zero.  Second, we assume that $f\in
C^\infty_0(\p M\times \R_+)$. Hence $u^f\in C^\infty(M\times \overline
\R_+)$ (see \cite{Lasiecka}), and $u^f(\widehat x, T)$ exists
pointwise.  Third, a function $v\in L^2(M)$ is interpreted as a
distribution $v\in \mathscr{D}'(M)$ by the formula
$$
  \bra v, \phi \cet = \int_M v \phi \,\mathrm{d}V, \quad \phi \in \mathscr{D}(M).
$$
{\newtext 
The delta distribution at $y\in M$ is
defined by $\bra \delta_y, \phi\cet = \phi(y)$ for $\phi\in \mathscr{D}(M)$. 
}

%
%


\section{Proofs} 
\label{proofSec} 
{\newtextt
We start with the proof of Proposition \ref{lemma_Control 1}.
The proof is a relatively direct consequence of Tataru's unique
continuation theorem and can be found e.g.\ in the case of Dirichlet
boundary conditions in  \cite[Lemma 2.1]{KL_edin}.
}
\begin{proof}[Proof of Proposition \ref{lemma_Control 1}]
{\newtextt Assume that a pair
\ba
(\psi,-\phi)\in (H^1(M)\times L^2(M))'=H^{-1}_0(M)\times L^2(M)
\ea
satisfy
the duality
\ba
\bra u^f(T),\psi\cet_{(H^1(M),H^{-1}_0(M))}+\bra u^f_t(T),-\phi\cet_{L^2(M)}=0
\ea
for all $f\in C^\infty_0(\partial M\times (0,T))$.
Note that $H^1(M)$ is the domain of the square root
of the operator $\A+cI$ when $c$ is large enough, denoted by
${\mathcal D}(\A^{1/2})$ and $H^{-1}_0(M)$ is the dual
$H^1(M)={\mathcal D}(\A^{-1/2})$.
Let
\ba
& & 
e_{tt}+\A e = 0\quad \hbox{ in }M\times (0,T),\\
& &\p_\nu e|_{\partial M\times (0,T)} = 0, \quad e|_{t=T} =\phi,\  e_t|_{t=T}  = \psi.
\ea
By \cite{Lasiecka}, 
 $e\in C^1([0,T],L^2(M))\cap C([0,T],H^{-1}_0(M))$ and
$e|_{\p M\times (0,T)}\in H^{-2/5-\e}(\p M\times (0,T)$, $\e>0$. Thus
 we have in sense of distributions 
\ba
0&=&\int_{M\times (0,T)}[u^f (e_{tt}+\A e)
-(u^f_{tt}+\A u^f) {e}]\,dV\,dt\,\\
&=&
\int_{M}(u_t^f(T)\,{\phi} -u^f(T)\, \psi)\,dV\,+
  \int_{M\times (0,T)} f\,e\,dS_x\,dt\,\\
&=&
\int_{M\times (0,T)} f\,  e\,dS_x\,dt
\ea
for all $f\in C^\infty_0(\partial M\times (0,T))$. This yields
that 
\ba
e|_{\partial M\times (0,T)}=
\partial_\nu e|_{\partial M\times (0,T)}=0.
\ea
To apply unique continuation for  
$e\in C([0,T],H^{-1}_0(M))$, let $\e>0$, and {\newtexttt let} $\eta\in 
C^\infty_0(\R)$ be a function supported on $(-1,1)\subset \R$ 
{\newtexttt whose} integral over $\R$ is one. Then
\ba
e_\e(x,t)=\int_\R e(x,t')\eta(\frac {t-t'}{\e})dt'
\ea
satisfies  
\ba
& & 
(\p_{t}^2+\A))e_\e = 0\quad \hbox{ in }M\times (\e,T-\e),\quad
\p_\nu e_\e|_{\partial M\times (\e,T-\e)} = 0
\ea
and
$e_\e\in C^\infty((\e,T-\e),H^{-1}_0(M))$. 
{\newtexttt By representing} $e_\e$ in terms of eigenfunctions of $\A$,
we see that $e_\e\in C^\infty((\e,T-\e),{\mathcal D}(\A^\infty))\subset
 C^\infty(M\times (\e,T-\e))$. Using
  Tataru's unique continuation theorem \cite{Ta1} we see that
if $0<\e<(T-\diam(M))/2$ then
$e_\e(T/2)=\p_te_\e(T/2)=0$. Hence $e_\e=0$ identically
on $M \times [0,T]$. When $\e\to 0$, we see that also $e$ vanishes identically 
and thus $\phi=\psi=0$.}
\end{proof}

\subsection{Proof of Theorem \ref{Main F}.}
\label{Sec: ProofOfFirstTheorem}
On $X$ we will study the minimization problem
\begin{eqnarray}
\label{eq: minimize}
  \min_{(h,a)\in X} \F(h,a,\alpha),
\end{eqnarray}
where $\alpha\in(0,1)$ {\newtexttt and} $\F$ is defined in equation \eqref{F intro
  definition}. {\newtextt By \cite{Lasiecka}, 
the map $h\mapsto u^h$ is continuous
$L^2(\p M\times [0,T])\to C([0,T];H^{5/6-\e}(M))$, $\e>0$. Thus  
$(h,a)\mapsto  \F(h,a,\alpha)$
is continuous map $X\to \R$.}

\begin{lemma}
\label{lemma_7} 
For any $\alpha\in(0,1)$ minimization problem \eqref{eq: minimize} has 
a unique minimizer $(h,\,a)\in X$. This minimizer is the unique solution to 
\beq
\label{Eq: main}
  (\alpha + L) \vc{h}{a} 
   = 
  \vc{PKf}{0 },
\eeq
where $L$ is defined in equation \eqref{Edef}.
Furthermore, $L\colon X\to X$ is non-negative, bounded, and
selfadjoint. 
\end{lemma}

\begin{proof}
We have 
\ba
\label{eq: minimize 3}
  \F(h,a,\alpha)&=& \bra f-Ph,K(f-Ph)\cet + \bra Ph-a,K(Ph-a)\cet \\ 
\nonumber
  & &\, + \, \bra\p_t a,K\p_ta\cet \\
  & &\, + \, \alpha(\bra h,h \cet+ \bra a,a \cet+ \bra \p_t a,\p_ta\cet),
\ea
Here $K$ and $P$ are selfadjoint in $L^2(M)$.
The dual $Y^*$ of the Hilbert space  $Y=H^{1}([0,2T];L^2(M))$
can be identified with $H^{-1}([0,2T];L^2(M))$. As
 $Q$ the inverse of $1-\p^2_t$, it can be considered
as an operator $Q\colon Y^*\to Y$. Thus using
$\p_t K\p_t\colon Y\to Y^*$, we can write 
\ba
\bra \p_t a,K\p_ta\cet_{L^2(M)} =-\bra \p_t K\p_ta, a\cet_{Y^*,Y}
=-\bra Q\p_t K\p_ta, a\cet_{Y}.
\ea
Thus, using the inner product on $X=L^2(M)\times Y$, we can rewrite $\F$ as
\begin{eqnarray}
\nonumber
  \F(h,a,\alpha) &=& \bra f, K f\cet  +2 \left\bra  \vc{h}{a}, 
\begin{pmatrix}
  1 & 0 \\
  0 & Q \\
\end{pmatrix} 
\vc{-PKf}{0}  \right\cet_{\!X}  \\
\label{eq: minimize 2}
& & \, + \, \left\bra  \vc{h}{a}, (\alpha  + L)   \vc{h}{a}  \right\cet_{\!X}\!.
\end{eqnarray}
As $Q\colon Y^\ast \to Y$ and $\p_t K \p_t \colon Y \to Y^\ast$ are
bounded, $L\colon X\to X$ is bounded. A direct calculation shows that $L$ is
self-adjoint.  Setting $f=0$ and $\alpha=0$ in equation \eqref{eq:
  minimize 2} shows that $L$ is non-negative.

Let us next observe  that $\F$ is strictly convex, so for the first claim 
it suffices to prove existence; by convexity, a local minimum is a global 
minimum, and by strict convexity, there is only one global minimum
(see \cite[ Theorem 38.C] {Zeidler:1984}). 
To prove existence, we recall that by  \cite{AMR}, $x\in X$ is a local minimum
of $\F\colon X\to \R$ provided that:
\begin{enumerate}
  \item[(i)] the first two Fr\'echet derivatives at $x$, $D_x\F(\xi)$ and $D_x^2 \F (\xi,\eta)$, exist and are continuous,
  \item[(ii)] $D_x\F=0$,  
  \item[(iii)] $u\mapsto D^2_x\F(\xi,\cdot)$ is a linear isomorphism $X\to X^\ast$, and 
  \item[(iv)] $D^2_x\F(\xi,\xi)>0$ for all $\xi\neq 0$.
\end{enumerate}
The sought Fr\'echet derivatives of $\F\colon X\to \R$ are 
\ba
 D \F_{h,a} ( \xi) &=&
2 \left\bra 
\begin{pmatrix}
  1 & 0 \\
  0 & Q \\
\end{pmatrix}   \vc{-PKf}{0 } + (\alpha + L)   \vc{h}{a}, \xi  \right\cet, \\
 D^2\F_{h,a}  (\xi,\eta) &=& 2\left\bra  \xi, (\alpha + L)  \eta  \right\cet, 
\quad \xi,\eta\in X. 
\ea 
It remains to prove that $\alpha + L$ is invertible, but if $(\alpha +
L)(\xi)=0$, then by non-negativity, $0\le \bra L\xi, x\cet_X = -\alpha \Vert
\xi \Vert^2_X \le 0$, and $\xi=0$.
\end{proof}


\begin{lemma}
\label{lemmaXYZ}
Iteration scheme \eqref{firstItScheme} converges to the unique
solution to equation \eqref{Eq: main}.
\end{lemma}

\begin{proof} 
Using $S$ and $\omega$ defined in definition \ref{cutDef}, we may rewrite
 equation \eqref{Eq: main} as
\begin{eqnarray*}
  (I-S)\vc{h}{a} = \frac{1}{\omega}\vc{PKf}{0},
\end{eqnarray*}
For a self-adjoint operator $B$, the norm satisfies
$\Vert B\Vert = \sup\{\vert \bra x,Bx\cet \vert : \Vert x\Vert = 1\}$.
Hence $\alpha + \bra L x, x\cet<\omega/2$, and by
non-negativity, 
\ba
\Vert S\Vert\le 1-\frac{\alpha}{\omega}<1.
\ea 
We may therefore iteratively solve $h,a$ by a convergent von 
Neumann sum.
\end{proof}

\begin{lemma} 
\label{finalLemma}
Minimizers $h(\alpha), a(\alpha) \in X$ for \eqref{eq: minimize}
satisfy
\begin{eqnarray*}
  \supp h(\alpha) &\subset&  B, \\
  a(\alpha) &\in & \operatorname{range} Q, \\
  \lim_{\alpha\to 0}
  \vc{  u^{a(\alpha)}(T)}{u^{a(\alpha)}_t(T)} &=& \vc{\chi_{N} u^f(T)}{0},
\end{eqnarray*}
where all limits are in $L^2(M)$. 
\end{lemma}

\begin{proof}
The first two claims follow by writing out \eqref{Eq: main}.
For the other results, let us define $Z\colon X\to \R$ by 
\begin{eqnarray*}
Z(h,a) &=& \frac{1}{2} \|\chi_N  u^f(T)-u^{Ph} (T) \|^2  + \frac{1}{4} \| u^{a}(T) - \chi_N u^f(T) \|^2\\
   & &\, + \, \|u^{a}_t(T) \|^2.
\end{eqnarray*}
To prove the last claim we show that for any $\varepsilon>0$ there exists 
an $\a(\varepsilon) \in (0,1)$ such that  
$Z(h(\alpha), a(\alpha)) < 4\varepsilon$ when $\alpha\in(0,\alpha(\varepsilon)).
$
{\newtext
Let us note that 
\begin{eqnarray*}
  u^{Pf}(T) = \chi_N u^{Pf}(T), \quad f\in L^2(\p M, [0,T]).
\end{eqnarray*}
}
Hence, for any $(h,a)\in X$,
\begin{eqnarray}
\nonumber
\F(Ph,a,\alpha) &=& \|(1-\chi_N) u^f(T) \|^2   + \|\chi_N u^f(T)-u^{Ph} (T) \|^2 \\
\nonumber
 & & +\, \| u^{Ph}(T) - u^a(T) \|^2 + \|u^a_t(T)\|^2 \\
\nonumber
 & & +\, \alpha(\|h\|^2+\|a\|^2+\|\p_ta\|^2).
\end{eqnarray}
It follows that for any $(h,a) \in X$ and $\alpha\in (0,1)$, 
\begin{eqnarray*}
  Z(Ph,a) &\le& \F(Ph,a,\alpha) - \|(1-\chi_N) u^f(T) \|^2.
\end{eqnarray*}
Here we have estimated the second term in $Z$ using the triangle
inequality and the inequality $(s+t)^2\le 2(s^2+ t^2)$.  Let us fix
$\varepsilon\in (0,1)$. By Proposition \ref{lemma_Control_2} there exists an
$h_\varepsilon\in L^2(B)$ such that
$$
  \|\chi_N u^{f}(T) - u^{Ph_\varepsilon}(T) \|^2 < \varepsilon,\\
$$
and by Proposition \ref{lemma_Control 1} there exists an 
$a_\varepsilon\in H^1(B)$ such that 
\ba
\| u^{a_\varepsilon}(T) - \chi_N u^{Ph_\varepsilon}(T) \|^2 &<& \varepsilon,\\
\| u_t^{a_\varepsilon}(T) \|^2 &<& \varepsilon.
\ea 
As $h_\varepsilon= Ph_\varepsilon$ we have
\ba 
  \F(h_\varepsilon,a_\varepsilon,\alpha) =  \|(1-\chi_N) u^f(T) \|^2    + 3\varepsilon + \alpha \left( \|h_\varepsilon\|^2+ \|a_\varepsilon\|^2 + \|\p_t a_\varepsilon\|^2 \right),
\ea 
and if $\alpha\in (0,\alpha(\varepsilon))$, where 
\ba 
  \alpha(\e) =\frac \varepsilon{1+\|h_\varepsilon\|^{2}+\|a_\varepsilon\|^2+\|\p_t a_\varepsilon\|^2
}, 
\ea
then the minimizer $h(\alpha),a(\alpha)$ of $\F$ satisfies
\ba 
  Z(h(\alpha), a(\alpha)) &\le & \F(Ph(\alpha),a(\a),\alpha) -\|(1-\chi_N) u^f(T) \|^2 \\
&\le & \F(h_\varepsilon,a_\varepsilon,\alpha) -\|(1-\chi_N) u^f(T) \|^2 \\
                            & < & 4\varepsilon. 
\ea 
\vspace{-1cm}
\end{proof}

\vspace{1cm}

\subsection{Proof of Theorem \ref{Cor 1}}
\label{focusWaveProof} 
{\newtextt 
Let us note that for any $\varepsilon>0$, then 
$\Gamma_k \subset B(\widehat z,\varepsilon)$ for sufficiently large $k$.
}
By Theorem \ref{Main F}, the following limits exist in $Y$,
\begin{eqnarray*}
   a(\alpha, \varepsilon) &=& \lim_{n\to \infty} a_n(\alpha, \varepsilon), \\
   a(\alpha, j, \varepsilon) &=& \lim_{n\to \infty} a_n(\alpha, j, \varepsilon).
\end{eqnarray*}
and
\begin{eqnarray}\label{eq: A1}
  \lim_{\alpha\to 0}
  \vc{  u^{a(\alpha, \varepsilon)}(T)}{u^{a(\alpha, \varepsilon)}_t(T)} &=& \vc{\chi_{N(\varepsilon)} u^f(T)}{0}, \\
\label{eq: A2}
  \lim_{\alpha\to 0}
  \vc{  u^{a(\alpha, j, \varepsilon)}(T)}{u^{a(\alpha,j,\varepsilon)}_t(T)} &=& \vc{\chi_{N(j,\varepsilon)} u^f(T)}{0},
\end{eqnarray}
where
\begin{eqnarray*}
  N(\varepsilon) &=& M(\p M,\widehat{T}-\varepsilon),  \\
  N(j,\varepsilon) &=& M(\p M, \widehat{T}-\varepsilon) \cup M(\Gamma_j, \widehat{T}+\varepsilon).
\end{eqnarray*}
We define 
$b(\alpha, j, \varepsilon) = \lim_{n\to \infty} b_n(\alpha, j,\varepsilon)$, 
whence
\begin{eqnarray*}
b(\alpha, j, \varepsilon) &=&  \varepsilon^{-\frac{m+1}{2}} \left(a(\alpha, j, \varepsilon) - a(\alpha, \varepsilon)\right).
\end{eqnarray*}

\begin{lemma} In $L^2(M)$,
$$
  \lim_{j\to \infty} \lim_{\alpha\to 0}
  \vc{  u^{b(\alpha,j,\varepsilon)}(T)}{u^{b(\alpha,j,\varepsilon)}_t(T)} =
\varepsilon^{-\frac{m+1}{2}}  \vc{ \chi_{J(\varepsilon)} u^f(T)}{0},
$$
where
$$
  J(\varepsilon) = M(\widehat z,\widehat{T}+\varepsilon) \setminus M(\p M, \widehat{T}-\varepsilon).
$$
\end{lemma}

\begin{proof} 
Since $a\mapsto u^a$ is linear, it suffices to prove that 
pointwise
$$
  \lim_{j\to \infty} \chi_{ M(\Gamma_j,\widehat{T} + \varepsilon) \setminus M(\p M, \widehat{T} - \varepsilon) } (x)  = \chi_{J(\varepsilon)}(x), \quad x\in M.
$$
This is clear for $x\in J(\varepsilon)$. If $x\notin J(\varepsilon)$
we claim that $x\notin M(\Gamma_j, \widehat{T}+\varepsilon)$ for
large $j$.  However, if $d(x,\widehat z)>\widehat{T}+\varepsilon$, then 
$$
  \Gamma_l \subset B\left(\widehat z, \frac{d(x,\widehat z)-\widehat{T}-\varepsilon}{2}\right)
$$
for large $l$.
For $y\in \Gamma_l$, $d(x,y) \ge d(x,\widehat z) - d(y,\widehat z) >
\widehat{T}+\varepsilon$, so $d(x,\Gamma_l) > \widehat{T}+\varepsilon$,
and $x\notin M(\Gamma_l, \widehat{T}+\varepsilon)$.
\end{proof}

The next Lemma show that $J(\varepsilon)$ are sets
that shrink onto $\widehat{x}$ in the case when $\widehat T<\tau(\widehat z)$.

\begin{lemma} [Properties of $J(\delta)$]
\label{Jdec}
 {\newtextt 
For any $\varepsilon>0$, there is a $\delta>0$ such that 
\begin{eqnarray}
\label{Jsubset}
 J(\delta) \subset B(\widehat{x}, \varepsilon).
\end{eqnarray}
Moreover, if $\widehat T<\tau(\widehat z)$}, then
$\{\widehat{x}\} \subset J(\delta)$ for all $\delta$
implying that $\bigcap_{\delta>0}J(\delta)=\{\widehat{x}\}$.
If $\widehat T>\tau(\widehat z)$ then 
$J(\delta)=\emptyset$ for $\delta$ small enough.
\end{lemma}

\begin{proof} Let us first prove (\ref{Jsubset}).
For a contradiction suppose that $\varepsilon>0$ and $x_1, x_2,
\ldots$ is a sequence such that
$$
  x_j \in J(1/j), \quad x_j \notin B(\widehat{x}, \varepsilon).
$$
As $M$ is compact, we can move onto a subsequence and assume that
$x_j$ converges to an $x \in M\setminus B(\widehat{x}, \varepsilon)$. Now
$d(x_j, \widehat z) \le \widehat T+1/j$ and $d(x_j, \p M) > \widehat T-1/j$, and as 
$x\mapsto d(x,\p M)$ is continuous,
$$
  d(x, \widehat z) \le \widehat T, \quad d(x, \p M) \ge \widehat T.
$$
Thus $\widehat T\le d(x,\p M) \le d(x,\widehat z) \le \widehat T$.
If  $\widehat T>\tau(\widehat z)$, we have $d(\widehat x,\p M)<\widehat T$
and obtain  a contradiction. Thus we can assume that
$\widehat T\leq \tau(\widehat z)$. Then the above inequalities yield that
$d(x,\widehat z)=d(x,\p M) = \widehat T$.
As $M$ is compact, there is a geodesic from $\widehat z$ to $x$ that realizes 
$d(\widehat x, \widehat z)$. Then $\eta$ also realizes $d(\widehat x, \p M)$, and $\eta$
must necessarily be normal to $\p M$ \cite{chavel}. Thus $\eta=\gamma$, and 
$\widehat x = x$; a contradiction. 
Thus (\ref{Jsubset}) is proven.

If  $\widehat T<\tau(\widehat z)$, then clearly 
$\widehat x\in J(\delta)$ for all $\delta>0$. 
On the other hand, if $\widehat T>\tau(\widehat z)$ then there is
$z'\in \p M$ such that $d(\widehat x,z')<d(\widehat x,\widehat z)$
and we see that $\widehat x\not \in J(\delta)$ for small $\delta$.
 Thus we have shown  that  $J(\delta)=\emptyset$ for $\delta$ small enough.
\end{proof}

\begin{proof}[Proof of Theorem \ref{Cor 1}]
{\newtextt Consider first the case when $\widehat T<\tau(\widehat z)$.
Then we observe that} the following limit exists
\begin{eqnarray}
\label{Cdefeqn}
  C(\widehat{x}) = \lim_{\varepsilon\to 0} \frac{\vol (J(\varepsilon))}{\varepsilon^{\frac{m+1}{2}}}
\end{eqnarray}
exists (see \cite{BKLS}). Here $\vol(A) = \int_A 1 \operatorname{d}V$ when
$A\subset M$. 
Let us also note that $B(\widehat{x}, \varepsilon/2)\subset
J(\varepsilon)$ so $\vol (J(\varepsilon))>0$. Thus, as
$u^f(T,\cdot)$ is continuous,
\begin{eqnarray*}
  \lim_{\varepsilon\to 0}   
\bra \frac{1}{\varepsilon^{\frac{m+1}{2}}} \chi_{J(\varepsilon)} u^f(T), \phi \cet 
  &=& C(\widehat x) \lim_{\varepsilon\to 0} \frac{1}{\vol (J(\varepsilon))} \int_{J(\varepsilon)} u^f(T,x) \phi(x)\mathrm{d}V(x) \\
  &=& \bra C(\widehat x) u^f(T,\widehat x) \delta_{\widehat x}, \phi\cet, \quad \phi\in \mathscr{D}(M).
\end{eqnarray*}
The result follows by \cite[Theorem 2.1.8]{HormanderI}. 

{\newtextt In the case when $\widehat T>\tau(\widehat z)$, 
$J(\e)=\emptyset$ for $\e$ small enough, and thus the 
limits (\ref{eq: A1}) and (\ref{eq: A2}) are the same. Hence
limit (\ref{eq: A0}) is zero.}
\end{proof}

\section*{Acknowledgements}
The research has been partially supported by Tekes project MASIT03,
EPSRC,
and the Academy of Finland Center of
Excellence programme 213476.

\end{document}